\documentclass[12pt,a4paper,psamsfonts]{amsart}
\usepackage{amssymb,amscd,amsxtra,calc}
\usepackage{cmmib57}

\theoremstyle{plain}
    \newtheorem{thm}{Theorem}[section]
    
    \newtheorem{claim}[thm]{Claim}
    \newtheorem{corollary}[thm]{Corollary}
    \newtheorem{lemma}[thm]{Lemma}
    \newtheorem{proposition}[thm]{Proposition}
    
    \newtheorem{theorem}[thm]{Theorem}

\theoremstyle{definition}

    \newtheorem{remark}[thm]{Remark}
\theoremstyle{remark}

    \newtheorem{setup}[thm]{}

\newcommand{\C}{\mathbb C}

\newcommand{\N}{\mathbb{N}}

\newcommand{\BPP}{\mathbb{P}}
\newcommand{\PP}{\mathbb P}
\newcommand{\Q}{\mathbb Q}

\newcommand{\R}{\mathbb R}

\newcommand{\Z}{\mathbb{Z}}

\newcommand{\Br}{\operatorname{Br}}

\newcommand{\diag}{\operatorname{diag}}

\newcommand{\id}{\operatorname{id}}

\newcommand{\NE}{\operatorname{NE}}
\newcommand{\Nef}{\operatorname{Nef}}
\newcommand{\NS}{\operatorname{NS}}

\newcommand{\Pic}{\operatorname{Pic}}
\newcommand{\Supp}{\operatorname{Supp}}
\newcommand{\Sing}{\operatorname{Sing}}

\newcommand{\Weil}{\operatorname{Weil}}

\newcommand{\ratmap}
{{\,\cdot\negmedspace\cdot\negmedspace\cdot\negmedspace\to\,}}

\newcommand{\Bs}{\mathrm{Bs}}

\begin{document}
\title
[Rationality of rationally connected threefolds]
{Rationality of rationally connected threefolds admitting non-isomorphic endomorphisms}
\author{De-Qi Zhang}
\address
{
\textsc{Department of Mathematics,
National University of Singapore} \endgraf
\textsc{10 Lower Kent Ridge Road, Singapore 119076, Singapore}
\endgraf
}
\email{matzdq@nus.edu.sg}
\begin{abstract}
We prove a structure theorem for non-isomorphic endomorphisms of
weak $\Q$-Fano threefolds (cf.~Theorem~\ref{wFano}), or more generally for threefolds with
big anti-canonical divisor.
Also provided is a criterion for a fibred rationally connected threefold to be rational.
As a consequence, we show (without using the classification)
that every smooth Fano threefold
having a non-isomorphic surjective endomorphism is rational.
\end{abstract}
\maketitle
\section{Introduction}
We work over the field $\C$ of complex numbers.
We shall prove Theorems \ref{Rat}, \ref{Fano} and \ref{wFano}.

\begin{theorem}\label{Rat}
Let $X$ be a $\Q$-factorial projective threefold
having only terminal singularities,
a surjective endomorphism $f: X \to X$ of degree $> 1$,
and a $K_X$-negative extremal contraction $X \to Y$
with $\dim Y \in \{1, 2\}$. Suppose either one of the following
three conditions:
\begin{itemize}
\item[(1)]
$X$ is Gorenstein.
\item[(2)]
$K_X^3 \ne 0$.
\item[(3)]
The ramification divisor $R_f$ is nonzero.
\end{itemize}
Suppose further that $X$ is rationally connected.
Then $X$ is rational.
\end{theorem}

\begin{theorem}\label{Fano}
Let $X$ be a smooth Fano threefold with a surjective endomorphism $f : X \to X$
of degree $> 1$. Then $X$ is rational.
\end{theorem}

Here is a sketch of the proof of Theorem \ref{Fano} using Theorem \ref{Rat}.
Let $X \to X_1$ be the composite of
blowdowns between smooth threefolds such that $X_1$ is a
primitive smooth Fano threefold in the sense of \cite{MM}.
The morphism $f$, replaced by some power, descends to a surjective morphism
$f_1 : X_1 \to X_1$, using the finiteness of extremal rays in the Mori
cone $\overline{\NE}(X)$.
If $\rho(X_1) = 1$, then $X_1 \cong \BPP^3$ by \cite{ARV} or \cite{HM}, done!
Suppose that $\rho(X_1) \ge 2$.
By \cite[Theorem 5]{MM}, $X_1$ has an extremal contraction of conic bundle type,
so $X_1$ is rational by Theorem \ref{Rat}.
Alternatively, as a referee remarked, such $X_1$ fits one of the four families (or else known to be rational;
cf. \cite{MM}) whose explicit description may be used to give another
probably simpler proof of Theorem \ref{Fano}.

In \cite{KX}, Koll\'ar-Xu constructed many examples of endomorphisms
$f : X \to X$ of degree $> 1$
on quotients $X$ of projective spaces, which are Fano varieties
of Picard number one with at worst terminal singularities;
some $X$ might be irrational by invoking David Saltman's famous
counter examples to Noether's problem.

As another type of example, one takes elliptic curve $E$ of period $\sqrt{-1}$
and the abelian variety $A = E^n$ with $n \ge 1$.
Let $\mu_4 \cong \langle \sqrt{-1} \rangle$
act diagonally on $A$.  Then $X = A/\mu_4$ has only log terminal singularities
and is rationally connected when $n \le 3$. For $m \ne 0, \pm 1$,
the multiplication map $m_A : A \to A$
($a \mapsto ma$) descends to an endomorphism $f : X \to X$ of degree $m^{2n} > 1$.
When $n = 3$, it is not known whether $X$ is rational.

These examples suggest that the Gorenstein requirement
might be necessary for Theorem \ref{Rat} or \ref{Fano}.

\par \vskip 1pc
For a uniruled threefold $X$ with $-K_X$ big and an endomorphism,
we have the following equivariant minimal model program
(MMP); see \ref{pFano} - \ref{rFano}
for applications.
As a referee remarked, the assertion(4) is well known and also follows from
\cite[Lemma 2.8]{PS}.

\begin{theorem}\label{wFano}
Let $X$ be a $\Q$-factorial projective threefold with only
terminal singularities, a big $-K_X$ and a surjective
endomorphism $f : X \to X$ of degree $q^3 > 1$.
Let $X = X_0 \ratmap X_1 \cdots \ratmap X_t$ be a composition
of $K_{X_i}$-negative flips and divisorial contractions
such that there is no $K_{X_t}$-negative extremal contraction of birational type.
Replacing $f$ by some power, we have {\rm (cf.~also~\ref{rFano}):}
\begin{itemize}
\item[(1)]
For every $1 \le j \le t$, the dominant rational map $f_j : X_j \ratmap X_j$
induced from $f$, is holomorphic. Set $f_0 := f$.
\item[(2)]
There is a $K_{X_t}$-negative extremal contraction $X_t \to Y$ with $\dim Y \le 2$.
The morphism $f_t$ descends to a surjective endomorphism
$h : Y \to Y$.
\item[(3)]
Either $\dim Y \in \{1, 2\}$, or
$X_t$ is a $\Q$-Fano threefold of Picard number one $($a consequence of MMP$)$.
In the latter case, all $f_i \ (0 \le i \le t)$ are polarized {\rm (cf. \ref{conv}(1)),} and further
$$f_i^* \, | \, N^1(X_i) = q \, \id.$$
\item[(4)]
Suppose further that $(X, \Theta)$ is klt weak $\Q$-Fano for some $\Theta \ge 0$.
Then for all $0 \le i \le t$, there are some $\Theta_i \ge 0$ such that
the pairs $(X_i, \Theta_i)$ are klt $\Q$-Fano.
\end{itemize}
\end{theorem}

\begin{remark} \label{rQ1}
Theorem \ref{Fano} was proved in \cite{uniruled} when the surjective
morphism $f: X \to X$ is assumed further to be polarized
(cf. \ref{conv}(1) below for its definition).
The author does not know whether Theorem \ref{Fano} follows from the result in \cite{uniruled},
but being polarized is a strong condition (nevertheless, see also
Theorem \ref{wFano}(3)). For instance, take any endomorphism
$g : Y \to Y$, then $f = g \times h : Y \times \PP^1 \to Y \times \PP^1$
is not polarized if $(\deg h)^{\dim Y} \ne \deg g$
(cf. \cite[Lemma 2.1]{nz2}).
\end{remark}

For an arbitrary uniruled variety $X$ with an endomorphism $f$, we like to
have an equivariant MMP: $X \ratmap X_m$ to reduce to the Fano fibration case.
To do so, we need to prove that the $K$-negative extremal rays appearing
in the composition $X \ratmap X_m$ of birational contractions are stabilized by $f$ and its
descents. This is not easy to prove because the Mori cone $\overline{\NE}(X)$ may have infinitely many
extremal rays and $X \ratmap X_m$ may involve flips and
hence may not be holomorphic.
Fortunately, for
Theorem \ref{wFano}, these difficulties are overcome  by Theorem \ref{fixext}.

We refer to
\cite{Fa}, \cite{Fav},
\cite{Mc}, \cite{ENS},
\cite{NZ}, \cite{nz2}, \cite{ICCM}, and \cite{uniruled}
for the recent development in the study of the dynamics of endomorphisms of
complex varieties.

\par \vskip 0.5pc \noindent
{\bf Acknowledgement.}
I would like to thank Hiromichi Takagi
for suggesting Lemma \ref{ne-im}(2),
Shigefumi Mori for examples of conic bundles with isolated
discriminant (cf. \cite[(1.1.1)]{MP}),
the audience of Kinosaki Symposium (October 2008)
and Pacific Rim Mathematical Association Conference (July 2009)
for the comments, and the referees for valuable suggestions
which make the presentation better.
This project is supported by an Academic Research Fund of NUS.

\section{Preliminary results}

\begin{setup}\label{conv}
{\bf Conventions} {\rm are as in Hartshorne's book, \cite{KMM} and \cite{KM}. Some more:}
\end{setup}

\par \noindent (1)
Every endomorphism $f$ of a projective variety in this paper is assumed to be surjective;
so it is finite by the projection formula.

\par
An endomorphism $f: X \to X$ is {\it polarized} if
there is an ample Cartier integral divisor $H$ such that $f^*H \sim qH$ for some
$q > 0$. So $\deg(f) = q^{\dim X}$.

\par \noindent (2)
Let $f : X \to X$ be an endomorphism and $\sigma_V : V \to X$
and $\sigma_Y : X \to Y$ morphisms. We say that
$f$ {\it lifts} to an endomorphism $f_V : V \to V$
if $f \circ \sigma_V = \sigma_V \circ f_V$;
$f$ {\it descends} to an endomorphism $f_Y : Y \to Y$
if $\sigma_Y \circ f = f_Y \circ \sigma_Y$.

\par \noindent (3)
Every boundary divisor in a pair is assumed to be
an effective $\Q$-divisor.

\par \noindent (4)
A pair $(X, \Theta)$
is called a {\it klt weak $\Q$-Fano variety}
(resp. {\it klt $\Q$-Fano variety})
if $X$ is a $\Q$-factorial normal projective variety,
the pair $(X, \Theta)$ has at worst klt singularities
as in \cite[Definition 2.34]{KM}
and $-(K_X + \Theta)$ is nef and big (resp. ample).
When $\Theta = 0$, the pair $(X, \Theta)$ is simply denoted as $X$.
Thus $X$ is a {\it klt weak $\Q$-Fano variety, or a klt
$\Q$-Fano variety} if so is $(X, 0)$.
A klt weak $\Q$-Fano variety
(resp. a klt $\Q$-Fano variety) with $(X, \Theta)$ terminal
(strengthened condition)
is called a {\it weak $\Q$-Fano variety}
(resp. a {\it $\Q$-Fano variety}).

\par \noindent (5)
Suppose that $X$ is a normal projective variety and $(X, \Theta)$ is terminal
or klt.
An {\it extremal ray} in the Mori cone $\overline{\NE}(X)$, the nef cone $\Nef(X)$ or some other convex cones
is in the sense of locally polyhedral cone.
A $(K_X + \Theta)$-{\it negative extremal ray} in $\overline{\NE}(X)$ is an extremal ray
whose intersection with $K_X + \Theta$ is negative.
By a {\it $(K_X + \Theta)$-negative extremal contraction} $X \to Y$, we mean the contraction of a
$(K_X + \Theta)$-negative extremal ray in $\overline{\NE}(X)$.
In particular, the Picard number $\rho(X) = \rho(Y) + 1$.
When $\Theta = 0$, such a contraction is simply called
an {\it extremal contraction}, or {\it divisorial contraction}, or {\it flip}, accordingly.

\par \noindent (6)
By a {\it conic bundle} $X \to Y$, we mean an extremal contraction
of relative dimension 1, where $X$ is normal projective with only terminal singularities
(so a general fibre is $\BPP^1$).

\par \noindent (7)
For a surjective morphism $f : X \to Y$,
the {\it discriminant} of $f$ is defined as $D(X/Y) := \{y \in Y \, ; \,
f$ is not smooth over $y\}$.
Let $D_s(X/Y)$ be the $s$-dimensional part of $D(X/Y)$, which
is also regarded as a reduced scheme.

\par \noindent (8)
For an endomorphism $\varphi : V \to V$ of a finite-dimensional
real vector space $V$,
$\rho(\varphi) := \max\{|\lambda| \, ; \, \lambda \in \C$
is an eigenvalue of $\varphi \}$ is the {\it spectral radius}.

\par \noindent (9)
For a $\Q$-factorial normal projective variety $X$,
denote $S(X) := \{$prime divisor $M \, ; \, M_{|M}$
is not pseudo-effective$\}$.
When $\dim X = 2$, our $S(X)$ is the set of negative curves on $X$.

\par \noindent (10)
For a normal projective surface $X$, denote by
$\Weil(X)$ the finite-dimensional $\R$-vector space of Weil divisors modulo
Mumford numerical equivalence: two Weil divisors are
{\it Mumford numerically equivalent} if their Mumford pullbacks
to a resolution are numerically equivalent.

\par \vskip 0.5pc
The (1) below is from \cite[Theorem 1.2.7]{MP}. For the (2), see
\cite[Theorem(3.5)]{Mo} (for the smooth case) and \cite[Corollary 2.4.2]{MP}
(originally due to Cutkosky).

\begin{lemma}\label{smoothbase}
Let $X \to Y$ be a conic bundle with $(\dim X, \dim Y)$
$= (3, 2)$. Then we have:
\begin{itemize}
\item[(1)]
$Y$ has at worst Du Val singularities of type $A_n$.
\item[(2)]
If $X$ is Gorenstein, then $Y$ is smooth.
\end{itemize}
\end{lemma}

\par
The following are known results about the rationality of threefolds;
see \cite[\S 2.2 - 2.3]{Isk}, \cite[\S 4.7]{Mi} and the references therein.

\begin{theorem}\label{Isk}
Let $X$ be a $\Q$-factorial projective threefold with at worst
terminal singularities. Then $X$ is rational if either one of the following
three conditions is satisfied.
\begin{itemize}
\item[(1)]
There is an extremal contraction $X \to B \cong \BPP^1$
such that a general fibre $F$ is a del Pezzo surface
with $K_F^2 \ge 5$.
\item[(2)]
$X$ is smooth. There is an extremal contraction $X \to Y$
onto a rational surface $Y$,
which is a non-standard conic bundle {\rm (cf. $\ref{rIsk}$ below).}
\item[(3)]
$X$ is smooth. There is an extremal contraction $X \to Y$
onto a rational surface $Y$,
which is a standard conic bundle. Either $Y = \BPP^2$
and the discriminant $D : = D(X/Y)$ has $\deg(D) \le 4$, or
there is a $\BPP^1$-fibration on $Y$ such that $F . D \le 3$ for
a general fibre $F$.
\end{itemize}
\end{theorem}

\begin{remark}\label{rIsk}
\begin{itemize}
\item[(i)]
In (2) and (3) above, $Y$ is smooth (cf.~Lemma \ref{smoothbase}).
According to \cite[4.7]{Mi},
a conic bundle $\pi: X \to Y$ is {\it standard}
if $\Pic(X) = \pi^* \Pic(Y) \oplus \Z K_X$.

\item[(ii)]
In (2) above, by \cite[\S 4.7]{Mi}, `$X \to Y$ is non-standard' if and only if
`$X \to Y$ is a $\BPP^1$-bundle in the Zariski-topology';
such a $\BPP^1$-bundle is locally trivial and hence $X$ is rational, due to the
triviality of the Brauer group $\Br(Y)$ for smooth rational surfaces.

\item[(iii)]
The (3) above was proved by Iskovskikh \cite[Theorem 1]{Isk-Duke};
see his survey \cite[\S 2.3, Theorem 8]{Isk} for its variation and references.
\end{itemize}
\end{remark}

For the first part of the result below, see \cite[Lemma 4.1, Remark 4.2]{Mi};
the last part follows from the first part and Riemann-Roch theorem
as proved in \cite[Lemma 2.3]{CCZ}.

\begin{lemma}\label{Dnc}
Let $X \to Y$ be a conic bundle with $X$ smooth and $(\dim X, \dim$ $Y)$
$= (3, 2)$. Then the discriminant $D = D(X/Y)$ is a divisor of normal
crossings, and every smooth rational component of $D \ ($if exists$)$ meets at least
two points of other components. In particular,
we have $|K_Y + D_i| \ne \emptyset$,
when $Y$ is rational and $D_i$ is
a connected component of $D$.
\end{lemma}

\begin{lemma}\label{newD}
Let $\pi: X \to Y$ be a conic bundle with $(\dim X, \dim Y) = (3, 2)$
and $Y$ a rational surface.
Then either the generic fibre is rational over $\C(Y)$
and hence $X$ is rational,
or there is a conic bundle $\pi' : X' \to Y'$ such that:
\begin{itemize}
\item[(1)]
$X' \, ($and hence $Y')$ are smooth, and there are
birational morphisms $\sigma_x: X' \to X$ and $\sigma_y : Y' \to Y$
satisfying $\pi \circ \sigma_x = \sigma_y \circ \pi'$.
\item[(2)]
No connected component $C_i'$ of $D(X'/Y')$ is contracted to a point by $\sigma_y$.
\item[(3)]
$($The $1$-dimensional part$)$ $D_1(X/Y) = \sigma_{y*} D(X'/Y') \ne \emptyset$.
\item[(4)]
$D_1(X/Y)$ and $D(X'/Y')$ have the equal numbers of connected components.
\end{itemize}
\end{lemma}

\begin{proof}
(1) is actually proved in \cite[Theorem 4.8]{Mi}.
If (2) is not true, $C_i'$ is contracted to a point
on $Y$ which is at worst of Du Val singularity of type $A_n$ (cf. Lemma \ref{smoothbase})
and hence $C_i'$ is a rational tree, which contradicts Lemma \ref{Dnc}.
The (3) is true
because every reducible fibre over some $d \in D(X/Y)$
should be underneath only reducible fibres over some $d' \in D(X'/Y')$ and
note that $\sigma_y : Y' \to Y$ is the blowup along $D(X/Y)$; see the
construction in \cite[Theorem 4.8]{Mi};
note also that $(\pi')^*E$ is irreducible for every prime divisor $E \subset Y'$
(and especially for those $E \subset D(X'/Y')$).
The (4) is from (2) and (3).
\end{proof}

\begin{proposition}\label{D<4}
Let $X \to Y$ be a conic bundle
with $Y$ a rational surface.
Suppose that either $D_1(X/Y) = \emptyset$, or there is a
smooth rational curve $F$ on $\widetilde{Y}$, with
$\sigma : \widetilde{Y} \to Y$ a minimal resolution, such that
$Bs|F| = \emptyset$ and $F . \sigma^*D_1(X/Y) \le 3$.
Then $X$ is rational.
\end{proposition}

\begin{proof}
We may assume that $D_1(X/Y) \ne \emptyset$ (cf. Lemma \ref{newD})
and $F$ is a general member in $|F|$.
Let $Y'' \to \widetilde{Y}$ be the blow up of $m = F^2$ points on $F$
away from $F \cap \sigma^{-1}D(X/Y)$.
Take a blowup $X'' \to X$ so that the composite
$X'' \to X \to Y \ratmap Y''$ extends
to a morphism with general fibre $\BPP^1$. By the proof of
\cite[Theorem 4.8]{Mi}, we have a lifting $\pi' : X' \to Y'$
of $\pi$ as in Lemma \ref{newD} so that $Y' \to Y$ factors as
$Y' \to Y'' \to \widetilde{Y} \to Y$.
Denote by $D' = D(X'/Y')$ and by $D''$, $\widetilde{D}$, $D$
its pushforwards (as cycles) on $Y''$, $\widetilde{Y}$ and $Y$. Note that
$D = D_1(X/Y)$ by Lemma \ref{newD}.

\begin{claim}\label{D<4c1}
$\sigma^*D \ge \widetilde{D}$.
\end{claim}

We prove Claim \ref{D<4c1}.
Write $\sigma^*D - \widetilde{D} = \sum_i e_i E_i$ with $E_i$
irreducible and $\sigma$-exceptional. Lemma \ref{smoothbase} implies that
each $E_i$ is a $(-2)$-curve. By Zariski's lemma, we have only to show the assertion
that $0 \ge E_j . (\sigma^*D - \widetilde{D}) = -E_j . \widetilde{D}$
for every $E_j$. We need to consider the case where
$E_j$ is a component of $\widetilde{D}$.
Applying Lemma \ref{Dnc} to $D'$ and using Lemma \ref{newD},
we have $E_j . \widetilde{D} \ge 0$.
This proves Claim \ref{D<4c1}.

\par \vskip 0.5pc
We continue the proof of Proposition \ref{D<4}.
The proper transform $F'' \subset Y''$ of $F$ is a fibre
of a $\BPP^1$-fibration. Let $F' \subset Y'$ be the total transform of $F''$.
Then $F' . D' = F'' . D'' = F . \widetilde{D} \le F . \sigma^*D \le 3$.
Now Proposition \ref{D<4} follows from Theorem \ref{Isk}.
\end{proof}

\begin{lemma}\label{Dstab}
Let $X \to Y$ be a conic bundle with $\dim Y = 2$, and
$f : X \to X$ a surjective endomorphism descending to an endomorphism $h : Y \to Y$.
Then $h^{-1}(D) = D$ for the $(1$-dimensional part$)$ $D := D_1(X/Y)$,
so after replacing $f$ by some power, we have
$h^{-1}(D_i) = D_i$ for every irreducible component
$D_i$ of $D$.
\end{lemma}

\begin{proof}
We have $h^{-1}(D) \subseteq D$ since the reducibility of a fibre $X_d$
over $d \in D$ implies that of $X_{d'}$ for every $d' \in h^{-1}(d)$.
So $D \supseteq h^{-1}(D) \supseteq h^{-2}(D) \supseteq \cdots$.
Considering the number of irreducible components, we have $h^{-s}(D) =
h^{-s-1}(D)$ for some $s> 1$. Since $h$ is surjective and applying $h^s$
and $h^{s+1}$ we get $h^{\pm}(D) = D$. This proves the lemma.\
\end{proof}

The result below can be proved as in \cite[Lemma 8 and Proposition 10]{Ny02}.
The assumption on $X_i$ can be weakened to
being normal Moishezon by using the Mumford
intersection on surfaces.

\begin{lemma}\label{S}
Let $\varphi: X_1 \to X_2$ be a finite morphism between two
$\Q$-factorial normal
projective surfaces with the same Picard number $\rho(X_1) = \rho(X_2)$.
Then we have:
\begin{itemize}
\item[(1)]
$\varphi$ induces a bijection ${\varphi} : S(X_1) \to S(X_2)
\ (C \mapsto \varphi(C))$ with a well-defined inverse
mapping $C'$ to $\varphi^{-1}(C')$.
Further, $\varphi^*C_2 = \alpha C_1$,
whenever $C_2 \in S(X_2)$, $C_1 = \varphi^{-1}(C_2)$
and $\alpha = \sqrt{\deg(\varphi) \, C_2^2/C_1^2}$.
\item[(2)]
Suppose that $\deg(\varphi) > 1$. Then $|S(X_i)| < \infty$;
we may also assume that ${\varphi} | S(X_1) = \id$,
when $\varphi : X_1 \to X_2 = X_1$ is an endomorphism
and $\varphi$ is replaced by some power.
\end{itemize}
\end{lemma}

As a referee remarked, Lemma \ref{ne} and Lemma \ref{ne-im}(2) are well known
and the latter also follows from \cite[Lemma 2.8]{PS}. We give proofs for them for the convenience of readers.

\begin{lemma}\label{ne}
Let $(X, \Theta)$ be a klt weak $\Q$-Fano variety.
Then we have:
\begin{itemize}
\item[(1)]
$(X, \Theta_0)$ is
klt $\Q$-Fano for some $\Theta_0 \ge \Theta$.
\item[(2)]
The
cone $\overline{\NE}(X)$ has only finitely many
extremal rays $(K_X$ negative or non-negative$)$.
\item[(3)]
$X$ is rationally connected.
\end{itemize}
\end{lemma}

\begin{proof}
(3) is proved in \cite[Theorem 1]{Zq}.
For (1), by \cite[Proposition 2.61]{KM}, $-(K_X + \Theta) \sim_{\Q} A_k + \frac{1}{k}E$
for some ample $\Q$-divisor $A_k$ and effective $\Q$-divisor $E$.
Choose $k >>0$ so that $(X, \Theta + \frac{1}{k}E)$ has at worst klt singularities;
see \cite[Corollary 2.35]{KM}.
Now $-(K_X + \Theta + \frac{1}{k}E) \sim_{\Q} A_k$ is ample. (1) is proved.
(2) is true by the cone theorem as in \cite[Theorem 3.7]{KM}.
\end{proof}

\begin{lemma}\label{ne-im}
\begin{itemize}
\item[(1)]
Let $\varphi : V \to W$ be a generically finite surjective
morphism between projective varieties.
Suppose that the cone $\overline{\NE}(V)$ has only finitely many extremal rays
$\R_{\ge 0}[C_i] \ (K_V$ negative
or non-negative$)$. Then $\overline{\NE}(W)$ has only finitely many extremal rays,
each of which is generated by some $\varphi(C_i)$.
\item[(2)]
Suppose that $(V, \Theta)$ is a klt weak $\Q$-Fano variety.
Let $\sigma : V \ratmap W$ be a $K_V$-flip or a $K_V$-divisorial
contraction. Then
$(W, \Theta_W)$ is klt $\Q$-Fano for
some $\Theta_W \ge 0$.
Hence both cones $\overline{\NE}(V)$ and $\overline{\NE}(W)$ have only finitely many
extremal rays. $($We remark that if $V$ is terminal then so is $W$; {\rm cf.} \cite[Corollary 3.42]{KM}$)$.
\end{itemize}
\end{lemma}

\begin{proof}
(1) is true, since $\varphi$ induces a surjective additive
map $\varphi_* : \overline{\NE}(V) \to \overline{\NE}(W)$.

For (2), consider the case where $\sigma$ is a flip.
Let $\pi : V \to Y$ and $\pi^+ : W = V^+ \to Y$ be the flipping contractions.
Let $H$ be an ample $\Q$-divisor on $Y$. By  Lemma \ref{ne},
we may assume that $-(K_V + \Theta)$ is ample after replacing $\Theta$.
Thus $-m(K_V + \Theta) - \pi^*H$ is very ample for $m >> 0$
and hence linearly equivalent
to an irreducible divisor $G_m$
Then $(V, \Theta + \frac{1}{m} G_m)$ is klt
by \cite[Corollary 2.31]{KM}.
Let $\Theta_W' \subset W$ be the proper transform of
$\Theta + \frac{1}{m}G_m$.
Then $W$ is $\Q$-factorial and $(W, \Theta_W')$ is klt by \cite[Proposition 3.37,
Corollary 3.42]{KM}. Also $-(K_W + \Theta_W')$ is $\Q$-linearly equivalent to the nef
and big divisor $(\pi^+)^*H/m$.
Now (2) follows from  Lemma \ref{ne}.

Consider the case where $\sigma$ is divisorial.
By Lemma \ref{ne}, we may assume that
$(V, \Theta)$ is klt $\Q$-Fano.
Take a small ample $\Q$-divisor $H_W \subset W$ such that $-(K_V + \Theta) - \sigma^*H_W$ is still ample
and hence it is $\Q$-linearly equivalent to some $\Delta > 0$
with $(V, \Theta + \Delta)$ klt by the argument in the previous case.
Then $-(K_V + \Theta + \Delta) \sim_{\Q} \sigma^*H_W$.
Set $\Theta_W := \sigma_*(\Theta + \Delta)$. Now
$-(K_W + \Theta_W) \sim_{\Q} \sigma_*\sigma^*H_W = H_W$ which is ample,
and $K_V + \Theta + \Delta \sim_{\Q} \sigma^*(K_W + \Theta_W)$.
Thus $(W, \Theta_W)$ is klt (and $\Q$-Fano) because so is $(V, \Theta + \Delta)$.
This proves the lemma.
\end{proof}

Below is Birkhoff's generalization of the Perron-Frobenius theorem (cf. \cite{Bi}).

\begin{theorem}\label{PF}
Let $C$ be a strictly convex closed cone of a finite-dimensional
real vector space $V$ such that $C$ generates $V$ as a vector space.
Let $\varphi: V \to V$ be an endomorphism such that $\varphi(C) \subseteq C$.
Then the spectral radius $\rho(\varphi)$ is an eigenvalue
of $\varphi$ and corresponds to an eigenvector in $C \setminus \{0\}$.
\end{theorem}

\begin{lemma}\label{radius}
Let $V$ be a normal projective surface and $f : V \to V$ a surjective endomorphism.
Then $\deg(f) \le \rho(f^*)^2$, where $\rho(f^*)$ is the spectral radius
of $f^* | \Weil(V)$.
\end{lemma}

\begin{proof}
By Theorem \ref{PF}, there is a nonzero $\R$-Cartier nef divisor $L \in \Weil(V)$ such
that $f^*L \equiv \rho L$ where $\rho := \rho (f^*)$. Take a Jordan canonical basis
$\{L_i\}_{i \ge 1}$ for $f^*|\Weil(V)$. Write $f^*L_i = \lambda_i L_i + L_{i-1}$
(set $L_0 := 0$). Note that $|\lambda_i| \le \rho$.
We may assume that $L . L_i = 0$ for all $i < s$
and $L . L_s \ne 0$.
Now $\deg(f) L . L_s = f^*L . f^*L_s = \rho L . (\lambda_s L_s + L_{s-1})
= \rho \lambda_s L . L_s$. So $\deg(f) = \rho \lambda_s \le \rho^2$.
\end{proof}

\begin{lemma}\label{scalar}
Let $X = X_0 \ratmap X_1 \cdots \ratmap X_r$ be a composition
of flips and divisorial contractions corresponding to
$(K_{X_i} + \Theta_i)$-negative extremal rays on some $\Q$-factorial klt
pairs $(X_i, \Theta_i)$. Suppose that $f = f_0 : X \to X$ is a surjective endomorphism
and descends to surjective endomorphisms
$f_i : X_i \to X_i \ (1 \le i \le r)$.
Suppose further that for some
$0 \le k \le r$ we have
$(f_k^{s_k})^*|N^1(X_k) = q^{s_k} \, \id$ for some $q > 1$ and integer $s_k > 0$.
Then all $f_j$  $(0 \le j \le r)$ are polarized, and
$(f_j^{s})^*|N^1(X_j) = q^{s} \, \id$ for some integer $s > 0$.
\end{lemma}

\begin{proof}
Let $H' \subset X_k$ be an ample $\Q$-divisor and set $H := \sum_{i=0}^{s_k-1} (f_k^i)^* H'/q^i$.
Then $f_k^* H \equiv qH$. Thus $f_k$ is polarized by \cite[Lemma 2.2]{nz2}.

By the remark above, we only need to show the second assertion, and for that
we will prove by the ascending induction on $|k-i|$.
We may assume that $s_k = 1$, after replacing $f$ by some power.
It suffices to show that $f_i^*|N^1(X_i) = q \, \id$ for $i = k-1$ and $k+1$.
So we may assume that $(k, r) = (0, 1)$ or $(1, 1)$.

Case(1). $\sigma : X \to X_1$ is a divisorial contraction with $E$
the exceptional divisor. $E$ is irreducible by \cite[Proposition 5-1-6]{KMM},
and $f^*E = eE$ for some $e > 0$ because both $f$ and its descent $f_1$ are necessarily finite
morphisms.
If $k = 0$ so that $f_i^* | N^1(X_i) = q \, \id$ for $i = 0$ then
the same holds for $i = 1$ because $\sigma^* N^1(X_1)$ is a subspace
of $N^1(X)$ and $f_1$ is the descent of $f$.
Suppose that $k = 1$ so that $f_1^* | N^1(X_1) = q \, \id$.
Then $f^* | N^1(X) = \diag[e, q, \dots, q]$ with respect to a basis
consisting of $E$ and a basis of $\sigma^*N^1(X_1)$.
For an ample divisor $H' \subset X_1$, the divisor
$H := \sigma^*H'$ is nef and big and
$f^*H \equiv qH$, so $e = q$ by \cite[Lemma 2.1]{nz2}.
Thus $f^*|N^1(X) = q \, \id$. We are done.

Case(2). $\sigma: X \ratmap X_1$ is a flip. Then $\sigma$ induces an isomorphic linear map
between $N^1(X)$ and $N^1(X_1)$ compatible with the two actions
$f_i^*|N^1(X_i)$ with $i = 0, 1$; see the proof of \cite[Proposition 3.37]{KM}.
So the scalarity of one action implies that of the other.
This proves the lemma.
\end{proof}

The theorem below is the key in proving Theorem \ref{wFano}.

\begin{theorem}\label{fixext}
Let $X$ be a $\Q$-factorial projective threefold with only
terminal singularities, big $-K_X$ and a surjective endomorphism $f : X \to X$ of
degree $d > 1$. Let $\R_{\ge 0}[\ell]$ be a $K_X$-negative extremal ray.
Then, replacing $f$ by some power, we have:
\begin{itemize}
\item[(1)]
$f(\ell)$ is parallel to $\ell$, i.e., $f(\ell) \in \R_{> 0}[\ell]$.
\item[(2)]
Suppose that $\R_{\ge 0}[\ell]$ gives rise to a divisorial contraction
$\sigma : X \to X_1$, or a flip $\sigma: X \ratmap X^+ = : X_1$, or a Fano fibration
$X \to X_1$ with $\dim X_1 \le 2$.
Then $f$ descends to a surjective endomorphism of $X_1$.
Further, $-K_{X_1}$ is also big when $\sigma$ is  divisorial or a flip.
\end{itemize}
\end{theorem}

\begin{proof}
We only need to show the first assertion of Theorem \ref{fixext} because
the second follows from the first;
see \cite[Lemmas 2.12 and 3.6]{uniruled}.
Indeed, the bigness of $-K_X$ implies that of $-K_{X_1} = \sigma_*(-K_X)$
when $\sigma$ is divisorial; $-K_{X_1}$ is also big when $\sigma$
is a flip, because $\sigma$ is then isomorphic in codimension one.

We now prove Theorem \ref{fixext}(1).
We may assume that $\ell$ is an irreducible curve.
Write $-K_X = A + D$ with $A$ an ample $\Q$-divisor
and $D$ an effective $\Q$-divisor.
By \cite[Theorem 3.7]{KM}, for every $1 > \varepsilon > 0$,
the first assertion below is true.

\begin{claim}\label{fixextc1}
\begin{itemize}
\item[(1)]
The subcone $R_{\varepsilon} := \overline{\NE}(X)_{K_X + \varepsilon A < 0}$
of $\overline{\NE}(X)$ contains only finitely many extremal rays.
\item[(2)]
Suppose that $\R_{\ge 0}[\ell]$ gives rise to a divisorial contraction
$\sigma: X \to X_1$ with $E$ the exceptional $($necessarily prime$)$ divisor
and $\ell$ taken to be a fibre of $\sigma$.
Then the subspace $N^1(X) . E \subset N_1(X)$ has rank $\le 2$,
and the ray $\R_{\ge 0}[\ell]$ is extremal in $($the closure of$)$ the cone $\Nef(X) . E \subset N_1(X)$.
\end{itemize}
\end{claim}

For Claim \ref{fixextc1}(2), if $\sigma(E)$ is a point then
$N^1(X) . E = \R[\ell]$, done!
Suppose that $\sigma(E)$ is a curve. For any divisor $M$ on $X$
we have $(M - aA) . \ell = 0$ for some $a \in \R$ and hence
$M - aA = \sigma^*L$ for some divisor $L$ on $Y_1$; see \cite[Lemma 3-2-5]{KMM}.
Then $M . E = aA . E + b \ell$ with $b = L . \sigma(E)$.
Thus $N^1(X) . E$ is spanned by $A . E$ and $\ell$.
Claim \ref{fixextc1} (2) is true because
the ray there, is extremal in the bigger cone $\overline{\NE}(X)$.
This proves Claim \ref{fixextc1}.

\vskip 0.5pc
We return to the proof of Theorem \ref{fixext}.

Note that every $\R_{\ge 0}[f^s(\ell)]$ ($s \in \Z$) is an extremal ray
in $\overline{\NE}(X)$ ($K_X$ negative or non-negative); see \cite[Lemma 2.11]{uniruled}.
If $f^s(\ell)$ is contained in $R_{\varepsilon}$ for infinitely many $s = s_1, s_2, \dots$,
then by Claim \ref{fixextc1}, we may assume that $f^{s_1}(\ell)$ is parallel
to $f^{s_2}(\ell)$ for some $s_2 > s_1$. Thus $f^{s_2-s_1}(\ell)$ is parallel to
$\ell$ because $f_*$ is an automorphism of the vector space $N_1(X)$. We are done.

Therefore, we may assume that $f^s(\ell)$ is not contained in $R_{\varepsilon}$ for
all $s \in \Z \setminus \{0\}$ after $f$ is replaced by some power.
Note that $\cup_{s > 0} \, f^s(\ell) \subset D$ because:
$$\begin{aligned}
-(K_X + \varepsilon A) &= (1 - \varepsilon)A + D, \\
0 \ge - f^s(\ell) . (K_X + \varepsilon A) &= (1 - \varepsilon) f^s(\ell) . A + f^s(\ell) . D
> f^s(\ell) . D.
\end{aligned}$$

Case (1). $R_{\ge 0}[\ell]$ gives rise to a Fano fibration $\sigma : X \to X_1$
with $\dim X_1 \le 2$ and $\ell$ chosen to be in a general fibre of $\sigma$.
Then $f(\ell)$ is not contained in $D$. This is a contradiction.

Case (2). $R_{\ge 0}[\ell]$ gives rise to a divisorial contraction $\sigma : X \to X_1$
with $E$ the exceptional (prime) divisor
and $\ell$ taken as a fibre of $\sigma$.
Since $f^s(\ell) . D < 0$ ($s \ne 0$), we have
$f(\ell) . D_1 < 0$ for an irreducible component $D_1$ of $D$.
Every fibre $\ell'$ of $\sigma$ is parallel to $\ell$,
so $f(\ell') . D_1 < 0$. Thus $f(\ell') \subset D_1$ and hence
$f(E) = D_1$. Now $f^2(\ell) . D < 0$ implies that $f^2(\ell) . D_j < 0$
for some irreducible component $D_j$ of $D$. As argued above, $f^2(E) = D_j$.
Since $D$ has only finitely many components, we have $f^{r_2}(E) = f^{r_1}(E)$ for
some $r_2 > r_1$. Thus for fibres $\ell_t$ of $\sigma$,
we have $f^{r_2}(\ell_t) = f^{r_1}(m_t)$ for some irreducible curve $m_t \subset E$.
So $f^r(\ell_t)$ (and hence $f^r(\ell)$) is parallel to $m_t$, where $r : = r_2 - r_1 > 0$.
For an extremal curve $\ell''$, denote
$$\Sigma_{\ell''} := \{C \subset X \, ; \, [C] \in \R_{\ge 0}[\ell'']\}, \,\,\,
U_{\ell''} := \cup_{\, C \in \Sigma_{\ell''}} \, C
$$
where these $C$ are irreducible curves. Then $E = U_{\ell}$.
By \cite[Lemma 2.11]{uniruled},
$$f^r(E) = f^r(U_{\ell}) = U_{f^r(\ell)} \supset \cup_t \,\, m_t \subset E.$$
The map $\{\ell_t\} \to \{m_t\}$ is a finite-to-finite map.
So $\{m_t\}$ is an infinite set and hence
$$f^r(E) = \overline{\cup_t \,\, m_t} = E.$$
By Claim \ref{fixextc1}, (the closure of) $\Nef(X) . E$ is generated by one or two
extremal rays (one being $\R_{\ge0}[\ell]$)
each of which is preserved by
$f^{2r}_*$. Thus $f^{2r}(\ell)$ is parallel to $\ell$. We are done.

Case (3). $\R_{\ge 0}[\ell]$ gives rise to a flip $X \ratmap X^+$
with flipping contractions $\pi : X \to Y_1$ and $\pi^+ : X^+ \to Y_1$.
We may assume that $f^{i}(\ell)$ and $f^{j}(\ell)$ are not
parallel for all $i \ne j$. In particular, $f^i(\ell)$'s are pair-wise
distinct. Replacing $f$ by some power, we may assume that
$$\cup_{s > 0} \,\, f^s(\ell) \, \subset \, \cup_{i=1}^r \, D_i \, \subseteq \, \Supp D$$
and $I_j := \{s > 0 \, | \, f^s(\ell) \subset D_j\}$ is an infinite set
for all $1 \le j \le r$. Thus
$$D_j = \overline{\cup_{s \in I_j} \,\, f^s(\ell)}.$$

\par
Each $\Sigma_{f^s(\ell)} = f^s(\Sigma_{\ell})$ is a union of
finitely many (and, when $s>>0$, the same number of) curves. Since
$f^{-s}(\Sigma_{f^s(\ell)}) = \Sigma_{\ell}$ (cf. \cite[Lemma 2.11]{uniruled}),
replacing $\ell$ by some $f^s(\ell)$,
we may assume that for $\ell_i := f^i(\ell)$ ($i \ge 0$)
we have $f^{-1}(\ell_{i+1}) = \ell_i$ and hence
$f^{-i}f^i(\ell_0) = \ell_0$. Now $f^{\pm}$ permutes $D_j$'s.
Indeed,
$$f(D_1) = f(\overline{\cup_{s \in I_1} \,\, f^s(\ell)})
= \overline{\cup_{s \in I_1} \,\, f^{s+1}(\ell)} \, \subset \, \cup_{i=1}^r \,\, D_i.$$
Hence $f(D_1) = D_{j_0}$ for some $j_0$. Also, by the choice of $j_0$ and $\ell_0$,
$$f^{-1}(D_{j_0})
= f^{-1}\overline{\cup_{s \in I_1} \,\, f^{s+1}(\ell)} =
\overline{\cup_{s \in I_1} \,\, f^{-1}f^{s+1}(\ell)} =
\overline{\cup_{s \in I_1} \,\, f^{s}(\ell)} = D_1.$$
Replacing $f$ by some power, we may assume that $f^{\pm}(D_j) = D_j$
for all $1 \le j \le r$, and also $f(\ell) . D_1 < 0$ after relabelling,
noting that only finitely many $f^s(\ell)$ are in $D \setminus \sum_{i=1}^r D_i$.
In particular, $f(\ell) \subset D_1$ and hence $D_1|D_1 \ne 0$.
Write $f^*D_1 = e_1 D_1$ for some integer $e_1 > 0$.
Then $f^*(D_1|D_1) = e_1 D_1|D_1$.

Consider the case $\deg(f|D_1 : D_1 \to D_1) = 1$. Then $e_1 = d > 1$.
We claim that $f|D_1$ induces an automorphism $h$ of positive entropy on the normalization
$\widetilde{D_1}$, the pullback of $-D_1|D_1$ is a
nef eigenvector of $h^*$ and $e_1$ equals the spectral radius $\rho (h^*)$ of
$h^* | \NS(\widetilde{D_1})$. Indeed,
applying Theorem \ref{PF} to the action of
$h^*$ on the nef cone of $\widetilde{D_1}$, we have
$h^*L = \rho (h^*) L$ for some nonzero $\R$-Cartier nef divisor $L$;
thus the pullback of $-D_1|D_1$ is a positive multiple of $L$ by
the uniqueness of eigenvalue of $h^*$ of modulus $> 1$ as in \cite[Theorem 3.2]{Mc}
applied to an equivariant desingularization of $D_1$. The claim is proved.

On the other hand,
$e_1$ (being
the spectral radius) must be a Salem number (and an algebraic integer
of degree $\ge 2$ over $\Q$)
and $e_1$ and $e_1^{-1}$
are roots of their common minimal polynomial (over $\Q$) dividing the
characteristic polynomial of $h^* | \NS(\widetilde{D_1})$;
see the proof of \cite[Theorem 3.2]{Mc} or \cite[Proposition 2.11]{-K}.
Indeed, $e_1 \in \N$ and hence $(x - e_1)$ is the minimal polynomial of $e_1$
over $\Q$, which is of degree one, absurd!

Therefore, $f|D_1 : D_1 \to D_1$ has degree $d_1 > 1$. Then $d = d_1e_1$.
Our $f|D_1$ induces a degree-$d_1$ endomorphism $h$ of the normalization $\widetilde{D_1}$
of $D_1$. After relabeling,
we may assume that no $\ell_i = f^i(\ell)$  ($i \ge 0$)
is included in the non-normal locus of $D_1$ or the image of the
ramification divisor $R_h$. Denote by $\widetilde{\ell_i} \subset \widetilde{D_1}$
the pullback of $\ell_i$. Then we have $h^* \widetilde{\ell_{j+1}} = \widetilde{\ell_j}$
for all $j \ge 0$. Let $N$ be the denominator of the intersection matrix
of the exceptional curves of a resolution of $D_1$. Then the Mumford intersection
satisfies
$C_1 . C_2 \in \frac{1}{N^2} \, \Z$ for all curves $C_i$ on $\widetilde{D_1}$.
Now
$$
d_1 \widetilde{\ell_i} . \widetilde{\ell_j} =
h^* \widetilde{\ell_i} . h^*\widetilde{\ell_j} =
\widetilde{\ell_{i-1}} . \widetilde{\ell_{j-1}}, \hskip 1pc
\frac{1}{N^2} \, \Z \, \ni \, \widetilde{\ell_i} . \widetilde{\ell_j} = \frac{1}{d_1^a}
\, \widetilde{\ell_{i-a}} . \widetilde{\ell_{j-a}}.
$$
Applying the above to $i = j = a >>0$ or $i = j-1 = a>>0$, we
get
$$\widetilde{\ell_i}^2 = \widetilde{\ell_{i+1}}^2 =
\widetilde{\ell_i} . \widetilde{\ell_{i+1}} = 0.$$
Thus $\widetilde{\ell_{i+1}}$ is parallel to $\widetilde{\ell_{i}}$,
by pulling back to a resolution of $D_1$
and applying the Hodge index theorem.
So $\ell_{i+1} = f^{i+1}(\ell)$ is parallel to $\ell_{i} = f^i(\ell)$.
Hence $f(\ell)$ is parallel to $\ell$. This proves Theorem \ref{fixext}.
\end{proof}

\begin{lemma}\label{etale}
Let $X = X_0$ and $X_1$ be $\Q$-factorial projective threefolds with only terminal
singularities and $f = f_0 : X \to X$ and $f_1 : X_1 \to X_1$ surjective endomorphisms
of the same degree $d > 1$.
Suppose that either $\sigma : X \to X_1$ is a $K_X$-divisorial contraction
or $\sigma: X \ratmap X_1$ is a $K_X$-flip, and that $f_1$ is the descent of $f$ in either case.
Then the ramification divisor $R_{f} = 0$,
if and only if $R_{f_1} = 0$,
if and only if $f_i$ is \'etale over $X_i \setminus \Sing X_i$ for both $i$.
\end{lemma}

\begin{proof}
Note that the second part follows from the first and the purity of branch loci.
If $\sigma: X \ratmap X_1$ is a flip, then $\sigma$ switches
$R_f$ and $R_{f_1}$ because $\sigma$ is an isomorphism in codimension $1$
and $f_1$ is the descent of $f$, so $R_f = 0$ if and only if $R_{f_1} = 0$.

Consider the case where $\sigma : X \to X_1$ is divisorial with
$E$ the exceptional (prime) divisor. Then $\sigma_*R_f = R_{f_1}$,
so $R_f = 0$ implies $R_{f_1} = 0$.

Suppose the contrary that $R_{f_1} = 0$ and $R_f \ne 0$.
Then $R_f = (e-1)E$, where $f^*E = eE$ with $e \ge 2$. Take a fibre $\ell$ of $\sigma$.
Then $\ell . K_X < 0$.
It is well known that $\ell . E < 0$.
Since $f$ and $f_1$ are compatible, $f^*\ell \equiv c \ell$ for some $c > 0$.
Now $d \ell . E = f^*\ell . f^*E = c e \ell . E$ implies that $c = d/e$.
Thus $f_* \ell = (d/c) \ell = e \ell$.
Now multiplying $\ell$ to the equality $K_X = f^*K_X + (e-1)E$ and by the projection
formula, we obtain
$$\ell . K_X = \ell . f^*K_X + (e-1) \ell . E
= e \ell . K_X + (e-1) \ell . E.$$
Thus $0 < -\ell . K_X = \ell . E < 0$, absurd.
So $R_f = 0$ if and only if $R_{f_1} = 0$.
\end{proof}

\section{Fibration-preserving endomorphisms}

In this section, we shall prove three lemmas used in the proof of the main theorems.

\begin{remark}
Every normal projective variety
$Y$ of dimension $\le 2$ and dominated by a rationally connected normal projective
variety $X$, is rational.
\end{remark}

\begin{lemma}\label{3-2 d_h=1}
Let $X$ be a
rationally connected projective threefold
with only terminal singularities,
and $\pi : X \to Y$ a conic bundle.
Suppose that $f : X \to X$ is an endomorphism of degree $d > 1$ and it descends to
an automorphism $h : Y \to Y$.
Then $X$ is rational.
\end{lemma}

\begin{proof}
By Proposition \ref{D<4}, it suffices to show the claim below.

\begin{claim}
The $(1$-dimensional part$)$ $D: = D_1(X/Y)$ is an empty set.
\end{claim}

We prove the claim.
Suppose the contrary that $D \ne \emptyset$.
Replacing $f$ by some power, we may assume, as in Lemma \ref{Dstab},
that $h^{-1}(D_i) = D_i$ for every
irreducible component $D_i$ of $D$. Thus $h^*D_i = D_i$ and $h^*D = D$ because
$\deg(h) = 1$.

Let $F$ be a very ample divisor on $Y$
away from $\Sigma := D_0(X/Y) \cup \Sing(D) \cup \pi(\Sing X) \cup \Sing(Y)$
and $h(\Sigma)$.
Set $s : = F . D > 0$.
Since the Picard number $\rho(X/Y) = 1$,
$X_F := \pi^*F$ is irreducible
and is a smooth ruled surface (admitting a ruling $X_F \to F$)
with $s$ singular fibres $S_i$ lying over the $s$ points $P_i$ in $D \cap F$.
Further, each $S_i$ consists of two intersecting lines $S_i(1) + S_i(2)$
(see \cite{Mo}, or \cite[\S 4.1, \S 4.8]{Mi}).

Note that the Picard number $\rho(X_F) = 2 + s$.
Denote $F' = h^*F$. We have $f^*X_F = X_{F'} := \pi^*F'$,
which is also a smooth ruled surface
with a ruling $X_{F'} \to F'$.
Our $f$ restricts to a finite morphism $f | X_{F'} : X_{F'} \to X_F$
of degree $d$ while the latter descends to the isomorphism
$h|F' : F' \to F$.

Now $D . F' = h^*D . h^*F = D . F = s$.
So $X_{F'}$ is a smooth ruled surface with the Picard number
$\rho(X_{F'}) = 2 + s = \rho(X_F)$. By Lemma \ref{S},
we have $(f|X_{F'})^*(S_i(j)) = d_{ij}S_i(j)'$ for some $d_{ij} \in \N$,
where $S_i(1)' + S_i(2)' = (f|X_{F'})^*S_i$ are the $s$ singular fibres of the ruling
$X_{F'} \to F'$ lying over the $s$ points $P_i' = D \cap F'$. Now
$-d = d S_i(j)^2 = ((f|X_{F'})^*S_i(j))^2 = (d_{ij}S_i(j)')^2 = -d_{ij}^2$,
so $d_{ij} = \sqrt{d} > 1$. Thus $h^*P_i = \sqrt{d} P_i'$ as 0-cycles,
with $P_i' = h^{-1}(P_i)$, contradicting $\deg(h) = 1$.
This proves the claim and Lemma \ref{3-2 d_h=1}.
\end{proof}

\begin{lemma}\label{3-2 d_h>1}
Let $X$ be a
rationally connected projective threefold
with only terminal singularities,
and $\pi : X \to Y$ a conic bundle.
Suppose that $f : X \to X$ is an endomorphism of degree $d > 1$ and it descends to
an endomorphism $h : Y \to Y$ of degree $e > 1$.
Then either $X$ is rational; or $d = e$, $X$ is non-Gorenstein,
$f$ is \'etale over $X \setminus \Sing X$, and $K_X^3 = 0$.
\end{lemma}

\begin{proof}
By Lemma \ref{smoothbase}, the
$Y$ is a rational surface with at worst Du Val singularities.
So $K_Y$ is Cartier and $K_{\widetilde Y} = \sigma^*K_Y$ if
$\sigma: \widetilde{Y} \to Y$ is a minimal resolution.

By Proposition \ref{D<4}, we may assume that the (1-dimensional part) $D : = D_1(X/Y) \ne \emptyset$.
Replacing $f$ by some power, we may assume that $h^{-1}(D_i) = D_i$
for every irreducible component $D_i$ of $D$, as in Lemma \ref{Dstab}. So $h^*D_i = d_iD_i$ for some $d_i \in \N$.
Thus $K_Y + D = h^*(K_Y + D) + G$
where $G$ is an effective divisor having no common component with $D$.
Inductively, for all $s > 0$, we have
$$K_Y + D = h^*(K_Y + D) + G =
(h^s)^*(K_Y+D) +
\sum_{i=0}^{s-1} (h^i)^*G.$$

In view of Proposition \ref{D<4}, we may impose an
extra assumption that $\sigma^*D . F' \ge 4$
where $F'$ is a general fibre of a $\BPP^1$-fibration on $\widetilde{Y}$ when
$\widetilde{Y} \ne \BPP^2$ and $F'$ is a line when $\widetilde{Y} = \BPP^2$.
We apply Lemma \ref{newD} and use the notation $X' \to Y'$ there.
We may assume that $Y' \to Y$ factors as $Y' \to \widetilde{Y} \to Y$.
Let $D' := D(X'/Y')$ and let $\widetilde{D} \subset \widetilde{Y}$
be the image of $D'$. By Lemma \ref{Dnc},
$|K_{Y'} + D'| \ne \emptyset$, so $|K_{\widetilde Y} + \widetilde{D}| \ne \emptyset$
and $|K_Y + D| \ne \emptyset$.

Take an ample divisor $H$ on $Y$.
Then for all $s > 0$,
$$H . (K_Y + D) = H . (h^s)^*(K_Y + D) + \sum_{i=0}^{s-1} H. (h^i)^*G \ge
\sum_{i=0}^{s-1} H. (h^i)^*G.$$
Thus $G = 0$ and $K_Y + D = h^*(K_Y + D)$.
If $K_Y + D \equiv 0$, then for the $F'$ above
we have $\sigma^*D . F' = \sigma^*(-K_Y) . F' = -K_{\widetilde Y} . F'\le 3$,
a contradiction to the extra assumption above.
So $K_Y + D$ is not numerically trivial.
By \cite[Theorem 7.1.1, cf. also Remark 7.1.2, Theorem 2.6.5]{ENS},
there is a fibration $Y \to B \cong \BPP^1$
(for $g(B) \le q(Y) = 0$)
with a general fibre $F$ such that
$$K_Y + D \sim_{\Q} aF$$
for some rational number $a > 0$. Since $h^*(K_Y+D) = K_Y + D$,
we have $h^*F \sim_{\Q} F$, so
$h$ descends to an automorphism
$h_B : B \to B$.
Now $0 = F. aF = F . (K_Y + D) \ge F . K_Y$.
Hence either $F \cong \BPP^1$ or $F$ is an elliptic curve.
If $F \cong \BPP^1$, then for $F' := \sigma^* F$, we have
$F' . \sigma^*D = F. (-K_Y) = 2$. This contradicts the extra assumption above.

We still have to consider the case where $F$ is elliptic
and $Y \to B$ is an elliptic fibration.

\begin{claim}\label{3-2c1}
\begin{itemize}
\item[(1)]
$h$ is \'etale over $Y \setminus \Sing Y$.
\item[(2)]
Every fibre of the elliptic fibration $Y \to B$ is irreducible.
\item[(3)]
We have $S(Y) = \emptyset$, so $\Sing Y \ne \emptyset$.
\item[(4)]
The Picard number $\rho(Y) = 2$.
\item[(5)]
If $\deg(f) = d = e = \deg(h)$, then $X$ is non-Gorenstein and
$f$ is \'etale over $X \setminus \Sing X$,
so $d K_X^3 = (f^*K_X)^3 = K_X^3$ and $K_X^3 = 0$.

\end{itemize}
\end{claim}

We now prove Claim \ref{3-2c1}.
Replacing $f$ by some power, we may assume that both $h$ and $h^{-1}$
stabilize every negative curve on $Y$ as in Lemma \ref{S}.
(1) is proved in \cite[Lemma 6.1.4]{ENS}.

(2) If $C$ is a curve in a reducible fibre then $C$ is a negative curve
and hence $h^*C = bC$ for some $b \in \N$, where
$b^2C^2 = (h^*C)^2 = eC^2$ and $b = \sqrt{e} > 1$.
Thus $h_B^*P = bP$ for the point $P$ over which lies the curve $C$.
This is impossible because $\deg(h_B) = 1$.

\par \vskip 0.5pc
Here is a note before continuing the proof of Claim \ref{3-2c1}:
Note that $0 = F . aF = F . (K_Y + D) = F . D$ and hence $D = \sum_{i=1}^s D_i$
with $D_i$ irreducible and the support of a fibre $F_i$.
So $D \equiv wF$ for some $w > 0$.
Since $K_{\widetilde Y}^2 = K_Y^2 = (aF - D)^2 = 0$,
the induced elliptic fibration $\widetilde{Y} \to \BPP^1$
(with a general fibre $\widetilde{F}$ say)
is relatively minimal, and $\sigma : \widetilde{Y} \to Y$ is just the
contraction of $(-2)$-curves in fibres of $\widetilde{Y} \to B$ into
Du Val singularities.

\par \vskip 0.5pc
(3) Note that $\widetilde{Y} = Y$ implies that there is a
$(-1)$-curve on $Y$
(noting that $\widetilde{Y}$ is rational and $K_{\widetilde Y}^2 = 0$).
So we have only to consider the case where there is a negative curve $E$
on $Y$. Then $h^*E = \sqrt{e}E$ as in (2).
Also $E$ is horizontal to the elliptic fibration $Y \to B$
by (2), so $E . F > 0$. Now
$e(E . F) = h^*E . h^*F = \sqrt{e} E . F$ and hence $e = \sqrt{e}$.
This is impossible for $e = \deg(h) > 1$ by the assumption.

(4) Since $\deg(h) = e > 1$ and $K_{\widetilde Y} = \sigma^*K_Y$ (and hence $K_Y$)
are not pseudo-effective, $\rho(Y) \ne 2$ would imply that
$h$ is polarized by some ample divisor $H$ so that
$h^*H \sim \sqrt{e} H$ as in \cite[Theorem 4.4.6]{ENS};
this leads to that $e (H . F) = h^*H . h^*F = \sqrt{e} H . F$, absurd.

(5) Note that $d/e = \deg(\ell \to f(\ell))$ for a general fibre
$\ell$ of $\pi : X \to Y$.
To distinguish, we write $f: X_1 = X \to X_2 = X$,
$h : Y_1 = Y \to Y_2 = Y$ and $\pi_i : X_i \to Y_i$
(identical to $\pi$) so that $\pi_2 \circ f = h \circ \pi_1$.

Suppose that $d = e$.
We assert then that $X_1$ is isomorphic to the
normalization $Z$ of $X_2 \times_{Y_2} Y_1$. Indeed,
$\pi_1 : X_1 \to Y_1$ and $f : X_1 \to X_2$ factor as
$X_1 \to Z \to Y_1$ and $X_1 \to Z \to X_2$,
where $Z \to Y_1$ and $Z \to X_2$
are the natural projections. Since $\pi_2 : X_2 \to Y_2$ has connected
fibres and $h$ is finite, $Z$ is irreducible. Note that the map $X_1 \to Z$ is finite
because so is $f: X_1 \to X_2$.
The map $X_1 \to Z$ is also birational
because $\deg(X_1/X_2) = \deg(Y_1/Y_2) = \deg(Z/X_2)$.
Thus the map $X_1 \to Z$ is an isomorphism and we can and will identify $X_1 = Z$.

By (1), $f : X \to X$ is \'etale in codimension $1$, so we have (*): $f$ is \'etale
over $X \setminus \Sing X$ by the purity of branch loci.
Thus $X$ is non-Gorenstein.
Indeed, a Gorenstein 3-dimensional terminal
singularity is an isolated hypersurface singularity and hence
has trivial local $\pi_1$ by a result of Milnor, so that
$f^s : X \to X$ is \'etale of degree $d^s$ for any $s > 0$ (by (*)),
contradicting the fact that the rationally connected variety $X$
has finite $\pi_1$.
This proves Claim \ref{3-2c1}.

\par \vskip 1pc
We continue the proof of Lemma \ref{3-2 d_h>1}.
Let $\rho$ be the spectral radius of $h^* | N^1(Y)$ and $L$
a nonzero $\R$-Cartier nef divisor such that $h^*L \equiv \rho L$.
Note that $\rho \ge \sqrt{e} > 1$ by Lemma \ref{radius}.
So $L$ is not parallel to $F$ and hence $L . F > 0$ by
the Hodge index theorem applied
to the pullbacks on $\widetilde{Y}$. Thus
$e (L . F) = h^*L . h^*F = \rho L . F$ and $\rho = e$.
Since $eL^2 = (h^*L)^2 = (eL)^2$, we have $L^2 = 0$.

\begin{claim}\label{3-2c2}
\begin{itemize}
\item[(1)]
We have $\overline{\NE}(Y) = \R_{\ge 0} L + \R_{\ge 0}F$
with $L$ and $F$ extremal rays.
\item[(2)]
$L$ $($replaced by its positive multiple$)$ is an integral Cartier divisor
with $h^0(Y, L) \ge 2$ and the Iitaka $D$-dimension $\kappa(Y, L) = 1$.
\item[(3)]
A multiple of $L$ is linearly equivalent to a general fibre of a
$\BPP^1$-fibration $\psi : Y \to C \cong \BPP^1$.
\end{itemize}
\end{claim}

We prove Claim \ref{3-2c2}.

(1) Note that $S(Y) = \emptyset$ and $L^2 = 0 = F^2$ imply that
$L$ and $F$ are extremal and hence (1) follows because the Picard number
$\rho(Y) = 2$.
Indeed, for $G = L$ or $F$,
decompose $G \equiv G_1 + G_2$ as sum of pseudo-effective divisors.
Let $G_i \equiv P_i + N_i$ be the Zariski-decomposition.
Since $S(Y) = \emptyset$, we have $N_i = 0$ and hence $G_i = P_i$
is nef. Now $0 = G^2 = G_1^2 + G_2^2 + 2G_1 . G_2$
implies that $G_i . G_j = 0$. Then, by the Hodge index theorem,
$G_1$ is parallel to $G_2$ and hence both $G_i$ are parallel to $G$.
So both $L$ and $F$ are extremal.

(2) We can choose $F$ and an ample divisor $H$ to be generators
of $N^1(Y)$. Replacing $L$ by its multiple, we may write
$L = \varepsilon F + uH$ with $\varepsilon = \pm 1$ and $u \in \R$.
Now $0 = L^2 = 2\varepsilon u F . H + u^2H^2$ implies
that $u = (-2\varepsilon F . H)/H^2 \in \Q$.
Thus we may assume that $L$ is integral and Cartier.

By the Riemann-Roch theorem and noting that $L . F > 0$ and
$\sigma^*K_Y = K_{\widetilde Y} = -v \sigma^*F$ for some $v > 0$
due to Kodaira's canonical divisor formula,
$$h^0(Y, L) = h^0(\widetilde{Y}, \sigma^*L) \ge \frac{L (L - K_Y)}{2} + 1 > 1.$$
So $\kappa(L) = 1$ because $L^2 = 0$.

(3) For the $L$ in (2), write $|L| = |M| + Fix$ with $Fix$ the fixed part.
Now $0 = L^2 \ge L . M \ge 0$ and  $0 = L^2 = L . M + L . Fix \ge M^2 + L . Fix \ge M^2 \ge 0$.
So $L . M = 0$ and hence $L$ is parallel to $M$ by the Hodge index theorem.
Also $M^2 = 0$ and hence $\Bs|M| = \emptyset$. Replacing $L$ by a multiple of $M$
we may assume that $L$ is a fibre of the fibration $\psi : Y \to C \cong \BPP^1$
derived from $\Phi_{|M|}$. Note that $K_Y = -vF$.
Since $\sigma^*L . K_{\widetilde Y} = L . K_Y = L . (-vF) < 0$,
our (3) is true.
This proves Claim \ref{3-2c2}.

\par \vskip 0.5pc
We continue the proof of Lemma \ref{3-2 d_h>1}.
By Claim \ref{3-2c2}(3), we may assume that $L$ is a general fibre of
$\psi$. Since the arithmetic genus $p_a(D_i) = 1$, we have $L . D > 0$.
Since $h^*L \equiv eL$, our $h : Y \to Y$ descends to an endomorphism
$h_C : C \to C$ of degree $e > 1$. By \cite[Theorem 5.1]{Fa}
and replacing $f$ by some power, we may assume that $h_C(c) = c$
for a general point $c \in C$. Our $h$ induces an endomorphism
$h|Y_{c} : Y_c \to Y_c$ of $Y_c := \psi^{-1}(c)$, and this endomorphism
is an automorphism for $\deg(h) = e = \deg(h_C)$.
Let $X_c = \pi^*Y_c$. Since $Y_c \equiv L$ by our choice of $L$,
we have $h^*Y_c \equiv eY_c$. So
$f : X \to X$ induces an endomorphism $f|X_c : X_c \to X_c$ of degree
$d/e$.

By Claim \ref{3-2c1}, we may assume that $d > e$.
Since $r := Y_c . D = L . wF > 0$, our $X_c$ is a smooth ruled surface
with a ruling $X_c \to Y_c \cong \BPP^1$ and $r$ singular fibres $S_i(1) + S_i(2)$
(two intersecting $(-1)$-curves) lying over the $s$ points $P_i \in Y_c \cap D$.
By Lemma \ref{S} and relabelling, we have $(f|X_c)^*S_i(1) = \sqrt{d/e} S_j(1)$.
Thus $(h|Y_c)^*P_i = \sqrt{d/e} P_j$ as $0$-cycles. This is impossible
for $\deg(h|Y_c) = 1$.
Lemma \ref{3-2 d_h>1} is proved.
\end{proof}

\begin{lemma}\label{3-1}
Let $X$ be a $\Q$-factorial rationally connected projective threefold
with only terminal singularities,
and $\pi : X \to Y$ an extremal contraction with $\dim Y = 1$.
Suppose that $f : X \to X$ is an endomorphism of degree $d > 1$.
Then either $X$ is rational; or $X$ is non-Gorenstein,
$f$ is \'etale over $X \setminus \Sing X$,
and $K_X^3 = 0$.
\end{lemma}

We now prove Lemma \ref{3-1}.
Since the Picard number $\rho(X) = \rho(Y) + 1 = 2$, we may assume that
$f$ preserves each of the two extremal rays of $\overline{\NE}(X)$, so
it descends to a surjective endomorphism $h : Y \to Y$
after replacing $f$ by its square.

\begin{claim}\label{3-1c1}
Suppose that $h$ is an automorphism.
Then $X$ is rational.
\end{claim}

\begin{proof}
Choose a general $y = y_0 \in Y$ such that $y_i := h^i(y_0)$ is not in $D(X/Y)$
or $\pi(\Sing X)$
for all $i \ge 0$. Set $X_i := \pi^*(y_i)$ which is a (smooth) del Pezzo surface.
Then $f^*X_{i+1} = X_i$ and
the restriction $f_i = f|X_i : X_i \to X_{i+1}$ is a finite surjective morphism of degree $d$.
By Theorem \ref{Isk}, we may assume that $K_{X_i}^2 = K_{X_0}^2 \le 4$ for all $i \ge 0$.
Note that the Picard number $\rho(X_i) = \rho(X_0)$.
We shall reach a contradiction late.

Let $N_i$ be the union of all negative curves (i.e., $(-1)$-curves) on $X_i$.
By Lemma \ref{S}, both $f$ and $f^{-1}$ induce natural bijections
(inverse to each other) between $N_i$ and $N_{i+1}$.
Indeed, $f_i^*N_{i+1}(j) = \sqrt{d} N_{i}(j)$
with $\deg(f_i) = d$,
if we label $N_i = \sum_j N_i(j)$ such that $f_i^{-1}N_{i+1}(j) = N_i(j)$.
Set $K_i = K_{X_i}$. Then
$K_{i} + N_{i} = f_i^*(K_{i+1} + N_{i+1}) + G_i$ where
$G_i$ is an effective divisor having no common components with $N_{i}$.
By iterating, for all $s > 0$, we have:
$$K_0 + N_0 = (f_s \circ f_{s-1} \circ \cdots f_0)^*(K_{s+1} + N_{s+1}) +
G_0 + \sum_{i=1}^s (f_{i-1} \circ \cdots f_0)^*G_i.$$

Since $K_i^2 \le 4$, our $N_i$ is a loop $N_i'$ plus a positive divisor $N_i''$,
so both $K_i + N_i'$ and $K_i + N_i$ are pseudo-effective (see \cite[Lemma 2.3]{CCZ});
indeed, a del Pezzo surface of degree $6$ contains a loop of $(-1)$-curves.
Multiplying the above displayed equality by an ample divisor
and letting $s \to \infty$, we see that $G_i = 0$ for almost all $i$.
Relabelling $X_i$, we may assume that $G_i = 0$ for all $i \ge 0$.
Thus, denoting $g_s := f_s \circ f_{s-1} \circ \cdots f_0$, we have
$$\begin{aligned}
K_0 + N_0 &= g_s^*(K_{s+1} + N_{s+1}) =
g_s^*(K_{s+1} + N_{s+1}') + g_s^*N_{s+1}'' \\
&= g_s^*(K_{s+1} + N_{s+1}') + d^{(s+1)/2}
\, g_s^{-1} N_{s+1}''.
\end{aligned}$$
Multiplying the equality by an ample divisor and letting $s \to \infty$,
we get a contradiction.
\end{proof}

By the above claim, we may assume that $\deg(h) > 1$.
Let $y_0 \in Y$ be a general $h$-periodic point which is not contained
in $\pi(\Sing X) \cup D(X/Y)$ or the branch locus of $h$ (see \cite[Theorem 5.1]{Fa}).
We may assume that $h(y_0) = y_0$ after replacing $f$ by some power.
Set $X_0 := \pi^*(y_0)$ which is a (smooth) del Pezzo surface.
Then the restriction $f_0 = f|X_0 : X_0 \to X_0$ is a finite morphism
of degree $d_0 := d/\deg(h)$.
If $d_0 > 1$, then $K_{X_0}^2 \ge 6$ by \cite{Ny02} or \cite[Theorem 3]{Zh02},
so $X$ is rational by Theorem \ref{Isk}.

Therefore, we may assume that $d_0 = 1$ and $\deg(h) = \deg(f) = d$.
So $f^*F \equiv dF$ for a general fibre $F$ of $\pi$.
To distinguish, we write $f: X_1 = X \to X_2 = X$,
$h : Y_1 = Y \to Y_2 = Y$ and $\pi_i : X_i \to Y_i$
(identical to $\pi$) so that $\pi_2 \circ f = h \circ \pi_1$.
As in the proof of Claim \ref{3-2c1}, $X_1$ is isomorphic to the
normalization $Z$ of $X_2 \times_{Y_2} Y_1$.

In particular, the ramification divisor $R_f$ is supported on fibres
and we can write $R_f \equiv bF$, noting that every fibre is irreducible
for $\rho(X/Y) = 1$.
If $b = 0$, i.e. $R_f = 0$, then
$K_X^3 = 0$ and $X$ is non-Gorenstein, as in the proof of Claim \ref{3-2c1}.
Thus we may assume that $b > 0$.

Since the Picard number $\rho(X) = 2$,
we have $\Nef(X) = \R_{\ge 0} L_1 + \R_{\ge 0} L_2$ with $L_i$ extremal rays.
We assert that a general fibre $F$ of $\pi$ is an extremal ray in $\Nef(X)$.
Indeed, if $F = G_1 + G_2$ is the decomposition
into nef divisors, then $0 = \ell . F = \ell . G_1 + \ell . G_2$
for a curve $\ell$ in $F$ and hence $\ell . G_i = 0$,
so $G_i \in \pi^*N^1(Y)$ and $G_i$ is parallel to $F$ (see \cite[Lemma 3-2-5]{KMM}).
The assertion is proved. So we may assume that $L_1 = F$ and hence $f^*L_1$ $\equiv dL_1$.
Also $L_1^2 = F^2 = 0$ (for later use).
We can also write $f^*L_2 \equiv \lambda L_2$ for some $\lambda > 0$.

Write $K_X \equiv a_1L_1 + a_2L_2$.
Since $-K_X|F \equiv -a_2L_2|F$ is ample we have $a_2 < 0$.
Now
$$K_X = f^*K_X + R_f \equiv (a_1dL_1 + a_2 \lambda L_2) + bL_1.$$
Comparing coefficients of $L_i$, we get
$$a_1 = \frac{b}{1-d} < 0, \hskip 2pc \lambda = 1.$$
So $-K_X$ is an interior point in $\Nef(X)$ and hence ample.
Also $f^*|N^1(X)$ can be diagonalized as $\diag[d, 1]$ with respect to
the basis $\{L_1, L_2\}$.
Write $\overline{\NE}(X) = \R_{\ge 0} \ell_1 + \R_{\ge 0} \ell_2 \subset N_1(X)$.
Since $N_1(X)$ is dual to $N^1(X)$, we have $B = \deg(f) (A^T)^{-1}$
if $A$ (resp. $B$) is the matrix representation of $f^*|N^1(X)$
(resp. $f^*|N_1(X)$) with respect to some basis.
Thus $f^*\ell_i \equiv r_i \ell_i$, where $\{r_1, r_2\} = \{1, d\}$.
We may assume that $\ell_1$ is contained in $F$ so that $\pi$ is the
contraction of the extremal ray $\ell_1$.
Then $\ell_1 . L_1 = 0$, so $\ell_1 . L_2 > 0$.
Thus $d (\ell_1 . L_2) = f^*\ell_1 . f^*L_2 = r_1 \ell_1 . L_2$, whence $r_1 = d$.
So $r_2 = 1$. Therefore, we have
$$f^*L_1 \equiv dL_1, \,\,\,\, f^*L_2 \equiv L_2, \,\,\,\,
f^*\ell_1 \equiv d \ell_1, \,\,\,\, f^*\ell_2 \equiv \ell_2.$$
$d (\ell_2 . L_2) = f^*\ell_2 . f^*L_2 = \ell_2 . L_2$ implies that
$\ell_2 . L_2 = 0$. Thus we have
$$\ell_i . L_i = 0 \,\, (i = 1, 2), \hskip 0.5pc \ell_i . L_j > 0 \,\, (i \ne j).$$

Let $\pi_2 : X \to Y_2$ be the contraction of
the (necessarily $K_X$-negative) extremal ray generated by $\ell_2$.
Then $f$ descends to an
endomorphism $h_2 : Y_2 \to Y_2$
and the Picard number $\rho(Y_2) = \rho(X) - 1 = 1$.
Hence $\dim Y_2 \ne 0$.
If $\dim Y_2 = 2$, the result follows from Lemmas \ref{3-2 d_h=1} and \ref{3-2 d_h>1}.

Consider the case $\dim Y_2 = 1$. Since $\ell_2 . L_2 = 0$,
the $L_2$ is parallel to a fibre of $\pi_2$.
Hence $L_2^2 = 0$. Take an ample $H$ and write
$H = u_1L_1 + u_2L_2$ with $u_i > 0$. Then $H^3 = 0$ for $L_i^2 = 0$,
absurd!

Consider the case $\dim Y_2 = 3$ so that $\pi_2$ is birational.
Let $H$ be an ample generator of $\Pic(Y_2)$. Then $\pi_2^*H$ is extremal in
$\Nef(X)$ as in the proof for the extremality of $L_1 = F$ above and since $\rho(Y_2) = \rho(X) - 1 = 1$.
Hence we may assume that $L_2 = \pi_2^*H$ because $\ell_2 . \pi_2^*H = 0$
and $\ell_2 . L_1 > 0$ imply that $\pi_2^*H$ is not parallel to $L_1$.
Set $h_2^*H = \alpha H$. Then $\alpha = 1$ because
$L_2 \equiv f^*L_2 = f^*\pi_2^* H = \pi_2^* h_2^* H = \alpha L_2$.
Thus $\deg(h_2) = (h_2^*H)^3/H^3 = 1$, whence $\deg(f) = 1$
for $h_2$ is the descent of $f$. This is a contradiction.
This proves Lemma \ref{3-1}.

\section{Proofs of Theorems and Corollaries}

In this section we shall prove Theorems \ref{Rat}, \ref{Fano} and \ref{wFano}
and their corollaries below.
In fact, Theorem \ref{Rat} is a consequence of Lemmas \ref{3-2 d_h=1}, \ref{3-2 d_h>1} and \ref{3-1}.

\begin{corollary}\label{pFano}
Let $X$ be a $\Q$-factorial projective threefold with only
terminal singularities and a surjective
endomorphism $f : X \to X$ of degree $q^3 > 1$.
Assume either one of the following three conditions:
\begin{itemize}
\item[(1)]
$X$ is Gorenstein.
\item[(2)]
$K_{X}^3 \ne 0$.
\item[(3)]
The ramification divisor $R_{f}$ is nonzero.
\end{itemize}
Assume further that $(*)$ : $X$ is rationally connected and $-K_X$ is big.
Then $X$ is rational, unless $(f^s)^* \, | \, N^1(X) = q^s \, \id$
for some integer $s \ge 1 \ ($and hence $f$ is polarized$)$.
\end{corollary}

The hypothesis $(*)$ in Corollary \ref{pFano} above is satisfied, if
$(X, \Theta)$ is klt weak $\Q$-Fano for some $\Theta$ (cf.~Lemma \ref{ne}).

\begin{corollary}\label{fFano}
With the situation as in Theorem $\ref{wFano} \, (3)$,
suppose that either $\dim Y = 0$ and $X_t$ is smooth, or
$\dim Y \in \{1, 2\}$ and at least one $X_i \ (0 \le i \le t)$
satisfies either one of the following three conditions:
\begin{itemize}
\item[(1)]
$X_i$ is Gorenstein.
\item[(2)]
$K_{X_i}^3 \ne 0$.
\item[(3)]
The ramification divisor $R_{f_i}$ is nonzero.
\end{itemize}
Suppose further that $X$ is rationally connected. Then $X$ is rational.
\end{corollary}

\begin{remark}\label{rCor}
\begin{itemize}
\item[(1)]
In the situation of Theorem \ref{wFano}, the ramification divisor $R_{f_k} = 0$
for some $k \in \{0, 1, \dots, t\}$, if and only if
$R_{f_i} = 0$ for all $i \in \{0, 1, \dots, t\}$,
if and only if $f_i$ is \'etale over $X_i \setminus \Sing X_i$; see Lemma \ref{etale}.
\item[(2)]
By the proof of Claim \ref{3-2c1}(5), in Theorem \ref{Rat} and Corollaries \ref{pFano} and \ref{fFano}, we have:
\newline
condition(1) $\Rightarrow$ condition(3); condition(2) $\Rightarrow$ condition(3).
\end{itemize}
\end{remark}

\begin{remark}\label{rFano}
Some remarks on Theorem \ref{wFano} related to the building blocks of endomorphisms:
\begin{itemize}
\item[(1)]
Every $X_i$ ($0 \le i \le t$) is $\Q$-factorial and has
only terminal singularities and big $-K_{X_i}$;
see \cite{KMM} or \cite[Proposition 3.37, Corollary 3.42]{KM}, and
Theorem \ref{fixext}(2).
\item[(2)]
Note that $-K_X$ is big and $X$ is rationally connected, whenever
$(X, \Theta)$ is klt weak $\Q$-Fano for some $\Theta$; see Lemma \ref{ne}.

\par \vskip 0.5pc \noindent
Now consider the case $\dim Y \in \{1, 2\}$ in Theorem \ref{wFano}(3).
\item[(3)]
Every general fibre of $X_t \to Y$ is $\BPP^1$ (when $\dim Y = 2$) or a smooth Fano surface
= del Pezzo surface (when $\dim Y = 1$), since $X_t$ has only terminal singularities.
\item[(4)]
Suppose further that $X$ is rationally connected.
If $\deg(h : Y \to Y) \ge 2$ and either the Picard number $\rho(Y) \ne 2$ or
$f_t^* | N^1(X_t) = q \, \id$
then $h$ is polarized (cf. \cite[Theorem 4.4.6]{ENS}
and \cite[Lemma 2.2]{uniruled}),
and hence the set of $h$-periodic points is dense in $Y$
by \cite[Theorem 5.1]{Fa}, so there is a surjective endomorphism
$h|X_y : X_y \to X_y$
for some general $y \in Y$, after
$h$ is replaced by some power. Here $X_y \subset X$ is the fibre
over $y$ and is $\BPP^1$ or del Pezzo.
\end{itemize}
\end{remark}

The remark below gives sufficient conditions to descend $f$ to a surjective
endomorphism $h : Y \to Y$ as in Theorems \ref{Rat} and \ref{wFano}.

\begin{remark} \label{rRat}
Let $\pi: X \to Y$ be a $K_X$-negative extremal contraction and $f : X \to X$ a
surjective endomorphism. Then $f$, replaced by some power, descends to
a surjective endomorphism $h : Y \to Y$ if any one of the following seven conditions
is satisfied.

\begin{itemize}
\item[(1)]
$X$ is a klt weak $\Q$-Fano variety.
\item[(2)]
$(X, \Theta)$ is a klt weak $\Q$-Fano variety for some $\Theta$.
\item[(3)]
The
cone $\overline{\NE}(X)$ has only finitely many extremal rays
($K_X$ negative or non-negative).
\item[(4)] $\dim Y \le 1$.
\item[(5)]
The Picard number $\rho(X) \le 2$.
\item[(6)]
$X$ is a $\Q$-factorial projective threefold with only terminal singularities
and big $-K_X$; the degree $\deg(f) \ge 2$.
\item[(7)]
$\dim Y < \dim X$; see \cite[Theorem 2.13, or Appendix]{uniruled}.
\end{itemize}
\end{remark}

\begin{setup}
{\bf Proof of Remark \ref{rRat}}
\end{setup}

The situation \ref{rRat}(6) is done by Theorem \ref{fixext}.
If $\dim Y \le 1$, then $\rho(X) = \rho(Y) + 1 \le 2$.
If $\rho(X) \le 2$, then $\overline{\NE}(X)$
has at most two extremal rays.
Thus by Lemma \ref{ne}, we have only to consider the situation
\ref{rRat}(3). Our $f^{\pm}$ permutes extremal rays
of $\overline{\NE}(X)$; see \cite[Lemma 2.11]{uniruled}. Hence $f^{\pm s}$ stabilizes every
such extremal ray for some $s > 1$. So Remark \ref{rRat} is true
(cf. \cite[Lemma 3.6]{uniruled} for the flip case).

\begin{setup}
{\bf Proof of Theorem \ref{Fano}}
\end{setup}

Since $X$ is Fano, it is rationally connected
and has only finitely many extremal rays $\R_{\ge 0} [C_i]$
all of which are $K_X$-negative; see Lemma \ref{ne}.
Let $X \to X_1$ be the composite of
blowdowns between smooth threefolds such that $X_1$ is a
primitive smooth Fano threefold in the sense of \cite{MM}.
Replacing $f$ by some power, we may assume that
$f$ descends to a surjective morphism $f_1 : X_1 \to X_1$.
This is because every smooth $X'$ appearing in between $X$ and $X_1$
is obtained by contracting a $K$-negative extremal ray,
the cone $\overline{\NE}(X')$ is generated by finitely many
extremal rays (the images of $\R_{\ge 0}[C_i]$) and a finite morphism
permutes these rays; see \cite[Lemma 2.11]{uniruled}.
If $\rho(X_1) = 1$, then $X_1$ is a smooth Fano threefold of Picard number one
having an endomorphism $f_1: X_1 \to X_1$ with $\deg(f_1) = \deg(f) > 1$,
so $X_1 \cong \BPP^3$ by \cite{ARV} or \cite{HM}, done!
If $\rho(X_1) \ge 2$,
by \cite[Theorem 5]{MM}, $X_1$ has an extremal contraction of conic bundle type,
so $X_1$ is rational by Theorem \ref{Rat}.
This proves Theorem \ref{Fano}.

\begin{setup}
{\bf Proof of Theorem \ref{wFano}}
\end{setup}

Let $X = X_0 \ratmap X_1 \cdots \ratmap X_t$ be a
composition of $K_{X_i}$ flips and divisorial contractions
so that there is no $K_{X_t}$-negative extremal contraction of birational type.

\begin{claim}\label{wFc1}
For every $0 \le i \le t$, the $X_i$ is
$\Q$-factorial uniruled with only terminal singularities,
$-K_{X_i}$ is big, $K_{X_i}$ is not nef, and there is a $K_{X_i}$-negative extremal contraction.
\end{claim}

We prove Claim \ref{wFc1}.
The bigness of $-K_X$ implies that of $-K_{X_i}$; see Theorem \ref{fixext}.
So we can write $-K_{X_i} = A_i + G_i$ with $A_i$ an ample $\Q$-divisor
and $G_i$ an effective $\Q$-divisor. Thus $X_i$ is uniruled
by \cite[Theorem 1]{MM}. For the other parts,
see \cite[Proposition 3.37, Corollary 3.42, Theorem 3.7]{KM}.
This proves Claim \ref{wFc1}.

\par \vskip 0.5pc
We return to the proof of Theorem \ref{wFano}.
Let $X_t \to X_{t+1} = Y$ be a $K_{X_t}$-negative extremal contraction.
Then $\dim Y \le 2$ by the choice of $X_t$.
Theorem \ref{wFano} (1) and (2) follow from Theorem \ref{fixext}.
For Theorem \ref{wFano} (3), if $\dim Y = 0$ then $\rho(X_t) = 1$
(so $f_t^* | N^1(X_t)$ is a scalar)
and hence $-K_{X_t}$ is ample because $K_{X_t}$ is not nef. Thus (3)
follows from Lemma \ref{scalar}.
(4) is proved in Lemma \ref{ne-im}.
This proves Theorem \ref{wFano}.

\begin{setup}
{\bf Proof of Corollary \ref{fFano}}
\end{setup}

Note that $\deg(f_t) = \deg(f) > 1$.
If $\dim Y = 0$, by the assumption, $X_t$ is a smooth Fano
threefold with a surjective endomorphism $f_t$ of degree $> 1$,
so $X_t \cong \BPP^3$ by \cite{ARV} or \cite{HM}, done!
Consider the case $\dim Y \in \{1, 2\}$.
By the proof of Theorem \ref{Rat} (in three lemmas), we may assume that $R_{f_t} = 0$.
Then $R_{f_i} = 0$ for all $i \in \{0, 1, \dots, t\}$ by Lemma \ref{etale}.
This proves Corollary \ref{fFano} (cf. Remark \ref{rCor}).

\begin{setup}
{\bf Proof of Corollary \ref{pFano}}
\end{setup}

In notation of Theorem \ref{wFano}, we may assume that $\dim Y \in \{1, 2\}$;
see also the first paragraph in the proof of Lemma \ref{scalar}.
Now
Corollary \ref{pFano} follows from Corollary \ref{fFano}.

\end{document}